\documentclass[twoside,reqno]{amsart}
\usepackage{amssymb}
\usepackage{epsfig}
\newtheorem{theorem}{Theorem}[section]

\newtheorem{definition}[theorem]{Definition}

\theoremstyle{remark}
\newtheorem*{remark}{Remark}

\newcommand{\refE}[1]{(\ref{E:#1})}
\newcommand{\refS}[1]{Section~\ref{S:#1}}
\newcommand{\refSS}[1]{Section~\ref{SS:#1}}
\newcommand{\refT}[1]{Theorem~\ref{T:#1}}

\newcommand{\refD}[1]{Definition~\ref{D:#1}}

\newcommand{\C}{\ensuremath{\mathbb{C}}}
\newcommand{\N}{\ensuremath{\mathbb{N}}}

\renewcommand{\O}{\ensuremath{\mathcal{O}}}
\renewcommand{\P}{\ensuremath{\mathbb{P}}}
\newcommand{\Z}{\ensuremath{\mathbb{Z}}}
\newcommand{\K}{\ensuremath{\mathbb{K}}}

\renewcommand{\i}{{\,\mathrm{i}\,}}
\newcommand{\im}{{\,\mathrm{Im}\,}}

\newcommand{\kno}{\K[X_0,X_1,\ldots,X_n]}

\newcommand{\Ls}{\mathrm{L}^2(M,L)}
\newcommand{\ghm}[1][m]{\Gamma_{hol}(M,L^{#1})}
\newcommand{\Tfm}[1][f]{T^{(m)}_{#1}}
\newcommand{\Qfm}[1][f]{Q^{(m)}_{#1}}
\newcommand{\gh}{\Gamma_{hol}(M,L)}
\newcommand{\Lm}[1][m]{L^{#1}}
\newcommand{\gul}{\Gamma_{\infty}(M,L)}

\renewcommand{\d}{\partial}
\newcommand{\db}{\overline{\partial}}
\newcommand{\qlb}{\ensuremath{(L,h,\nabla)}}
\newcommand{\w}{\ensuremath{\omega}}
\newcommand{\BT}{Berezin-Toeplitz}
\newcommand{\End}{\mathrm{End}}
\newcommand{\Cim}{C^{\infty}(M)}
\newcommand{\rank}{\mathrm{rank\;}}
\newcommand{\codim}{\mathrm{codim\,}}
\newcommand{\U}{\mathrm{U}}
\newcommand{\GL}{\mathrm{GL}}
\newcommand{\SL}{\mathrm{SL}}
\newcommand{\PGL}{\mathrm{PGL}}
\newcommand{\Con}{\mathbb{C}[X_0,\ldots,X_n]}
\newcommand{\Kon}{\mathbb{K}[X_0,\ldots,X_n]}
\newcommand{\Zero}{\mathcal{Z}}
\newcommand{\Ideal}{\mathcal{I}}

\begin{document}
\setcounter{page}{1}

\vspace*{-1cm}
\hspace*{\fill} Mannheimer Manuskripte 254

\hspace*{\fill} math/0005288
\vspace*{1cm}
\title[Singular projective varieties and quantization]
{Singular Projective Varieties and Quantization}
\author[M.~Schlichenmaier]{Martin Schlichenmaier}

\address{
        Department of Mathematics and 
           Computer Science\\
          University of Mannheim, D7, 27 \\
         D-68131 Mannheim \\
         Germany}

\email{schlichenmaier@math.uni-mannheim.de}

\begin{abstract}
By the quantization condition compact 
quantizable K\"ahler manifolds can be embedded into
projective space. In this way they become projective varieties.
The quantum Hilbert space of the Berezin-Toeplitz quantization 
(and of the geometric quantization) is the projective
coordinate ring of the embedded manifold.
This allows for generalization to  the case of 
singular varieties.
The set-up  is explained in the first part of the contribution. 
The second part of the contribution is of tutorial nature. 
Necessary notions, concepts, and results  of algebraic geometry 
appearing in this approach to quantization are explained.
In particular, the notions of projective varieties, 
embeddings, singularities, and
quotients appearing in geometric invariant theory 
are recalled.
\end{abstract}
\date{May 19, 2000}
\maketitle
\tableofcontents
\section*{Introduction}\label{S:intro}
Compact K\"ahler manifolds which are quantizable, i.e. which admit a
holomorphic line bundle with curvature form equal to the 
K\"ahler form (a so called quantum line bundle) are projective 
algebraic manifolds. This means that 
with the
help of the global holomorphic sections of a suitable tensor power of the 
quantum line bundle they can be embedded 
into a projective space of certain dimension.
Submanifolds of the projective space are always projective 
varieties, i.e. can be given as zero sets of finitely many 
homogeneous polynomials. 
As will be explained in this contribution the basic objects 
in the set-up of 
Berezin-Toeplitz (or equivalently geometric) quantization 
of quantizable compact K\"ahler manifolds can be completely 
described inside this algebraic-geometric context.
For example the quantum Hilbert space will be 
essentially the projective coordinate ring
of the variety.

By  definition  K\"ahler manifolds are nonsingular hence 
the varieties obtained are nonsingular.
But from the point of view of varieties the singular ones are on 
equal footing.
Hence one might expect that it is possible to find a direct way 
towards quantization of singular spaces by exploiting the 
theory of varieties.

In this contribution I do not present
a solution for the quantization of singular spaces. I will only 
explain the above mentioned path from compact quantizable
K\"ahler manifolds to projective varieties.
The quantization procedure I am considering is the Berezin-Toeplitz 
quantization, resp. the Berezin-Toeplitz deformation quantization.
This quantization procedure is adapted to the complex structure
which is a requirement for the fact that it can  be formulated in terms of 
complex algebraic geometry. I recall the results 
on this quantization  scheme in \refS{quant}.

The rest of the contribution is considered to be tutorial.
There is nothing new there, and everything 
is well-known for researchers working in algebraic geometry.
But I  hope that the collection of concepts and
results will be useful for researcher in quantization.
Some
concepts used elsewhere in this volume are explained.
In \refS{proj} and in \refS{sing} basic concepts of
algebraic geometry are introduced. 
First projective varieties are defined. 
Compactified moduli spaces are  candidates for
projective varieties.
The projective (homogeneous) coordinate ring
is discussed. It will turn out to be 
the quantum Hilbert space of the theory.
It incorporates the vector space of global holomorphic sections of
all tensor power of the quantum line bundle at once.
On this quantum Hilbert space the 
total Berezin-Toeplitz quantization 
operator operates. It is used to show 
that the quantization scheme has the correct semi-classical limit and
to prove the existence of 
an associated  deformation quantization.
As already pointed out, quantizable K\"ahler manifolds can be embedded 
with the help of the  very ample quantum line bundle into projective 
space (as complex manifolds, not necessarily as K\"ahler manifolds).
Such embeddings are discussed in detail in \refS{proj}.

Projective varieties are not necessarily smooth, they can have
singularities. After giving some examples of singular varieties in
\refS{proj} (e.g. the singular cubic curves) singularities are
treated in more detail in  \refS{sing}. Beside the definition of a 
singular point using the rank of the Jacobi matrix  of the defining
equations for the variety, a more intrinsic definition in
terms of the local ring $\O_{V,\alpha}$ of a point $\alpha$ on the
variety $V$ and the Zariski tangent space 
at the point $\alpha$ is given. 
In terms of algebraic properties of the local ring 
a hierarchy of types of singularities can be introduced.
As special examples  normal singularities are discussed.
Whereas on an arbitrary singular variety the set of singular
points can have codimension one, on a normal variety (i.e. on  a variety
where all local rings are normal rings, see the definition below)
this subset has codimension at least two.
The singularities of moduli spaces are very often normal singularities.

Typically, moduli spaces are  obtained by dividing out a group
action on a nonsingular variety.
The main question is whether it is possible  to define a geometric structure
on the orbit space,
i.e. whether there exists some algebraic geometric quotients.
The ``Geometric Invariant Theory (GIT)'' as developed by
Mumford \cite{GIT} gives a powerful tool how to deal with such
quotients.
If one considers only certain suitable subsets of points of the
variety the group is acting on (i.e. the subset of semi-stable, or
stable points) one obtains a good quotient
(which is also a categorical quotient), resp. a geometric quotient.
They will carry a compatible structure of a projective variety, resp. of 
an open 
subset of a projective variety.
This will be explained in \refS{quot}.
Also there the results on the relation with the symplectic quotients obtained
via moment maps and symplectic reduction due to  Kirwan, Kempf and Ness
will be explained. These results are taken from the appendix
to \cite{GIT}, written by  Kirwan.
Roughly speaking, the geometric quotient and the symplectic quotient
coincides on the regular points of the symplectic reduction
(see \refT{sympl} for a precise statement).
But in general the singularity structure will differ.

\section{From quantizable compact K\"ahler manifolds to
projective varieties}\label{S:quant}
Let $(M,\w)$ be a K\"ahler manifold, i.e. $M$ a complex manifold and
$\w$ a K\"ahler form on $M$. 
In this contribution I will only  consider
compact K\"ahler manifolds. If nothing else is said 
compactness is assumed.
A further data we need is the triple $\qlb$, with  a holomorphic
line bundle  $L$ on $M$, a hermitian metric $h$ on $L$ 
(with the convention that
it is conjugate linear in the first argument) and a connection
$\nabla$ compatible with the metric on $L$ and the complex structure.
 With respect to local holomorphic coordinates of the manifold
and with respect to a local holomorphic frame for the bundle
the metric $h$ can be given as 
\begin{equation}\label{E:locmet}
 h(s_1,s_2)(x)=\hat h(x)\overline{\hat s_1}(x)\hat s_2(x),
\end{equation}
where $\hat s_i$ is a local representing function for the
section $s_i$ ($i=1,2$) and $\hat h$ is a locally defined real-valued 
function on $M$.
The compatible connection is  uniquely defined and is given
in the local coordinates as
$\ \nabla=\d +(\d\log \hat h) +\db$.
The curvature of $L$ is defined as the two-form 
\begin{equation}\label{E:curvdef}
curv_{L,\nabla}(X,Y):=\nabla_X\nabla_Y-\nabla_Y\nabla_X-\nabla_{[X,Y]}\ ,
\end{equation}
where $X$ and $Y$ are vector fields on $M$.
In the local coordinates the curvature can be expressed as
$
curv_{L,\nabla}=\db\d\log \hat h=-\d\db\log \hat h
$.

A K\"ahler manifold $(M,\w)$ is called \emph{quantizable}
if there
exists  such a triple $\qlb$ which obeys
\begin{equation}\label{E:quantcond}
curv_{L,\nabla}(X,Y)=-\i\w(X,Y)\ .
\end{equation}
The condition \refE{quantcond} is called the (pre)quantum condition.
The bundle $\ \qlb\ $ is called a (pre)quantum line  bundle.
Usually we will drop $\nabla$ and sometimes also
$h$  in the notation.

For the following we assume  $(M,\w)$  to be a quantizable
K\"ahler manifold with quantum line bundle  $\qlb$.

There is an important observation.
If $M$ is a compact K\"ahler manifold which is quantizable then from
the prequantum condition \refE{quantcond}
 we obtain for  the Chern form of the line bundle
the relation
\begin{equation}
c(L):=\frac {\i}{2\pi}curv_{L,\nabla}=\frac {\w}{2\pi}\ .
\end{equation}
This implies that $L$ is a positive line bundle. In the terminology
of algebraic geometry it is an ample line bundle,
see \refD{ample} for the definition of ampleness.
By the Kodaira embedding theorem  $M$ can be embedded (as  algebraic
submanifold) into projective space $\P^N(\C)$ using
a basis of the global holomorphic
  sections $s_i$ of a suitable tensor power
$L^{m_0}$ of the bundle $L$
\begin{equation}
z\ \mapsto\  (s_0(z):s_1(z):\ldots:s_N(z))\in \P^N(\C)\ .
\end{equation}
These algebraic submanifolds  can be described
as  zero sets  of homogeneous polynomials,
i.e. they are projective varieties.
Note that the  dimension of the space $\ghm[m_0]$
consisting of the global holomorphic sections of $\Lm[m_0]$, can be determined
by the Theorem of Grothendieck-Hirzebruch-Riemann-Roch, see
\cite{GH}, \cite{SchlRS}.
By passing to the K\"ahler form $m_0\omega$ and to the 
associated  quantum 
line bundle $L^{m_0}$ we might assume that 
the sections of our 
quantum line bundle do already the
embedding (i.e that it is already very ample).

So even if we start with an arbitrary K\"ahler manifold  the
quantization condition will force the manifold to be an algebraic manifold
and we are in  the realm of algebraic geometry.
This should be compared with the fact that there are ``considerable more''
K\"ahler manifolds than algebraic manifolds.

In \refS{proj} I will explain what projective varieties are.
But first I like to introduce the quantum operator we are dealing with.
We take $\ \Omega=\frac 1{n!}\w^n\ $ as volume form  on $M$.
On the space of $C^\infty$ sections of the bundle $L$
 we have the scalar product
\begin{equation}
\langle\varphi,\psi\rangle:=\int_M h (\varphi,\psi)\;\Omega\  ,
\qquad
||\varphi||:=\sqrt{\langle \varphi,\varphi\rangle}\ .
\end{equation}
Let $\Ls$ be the  L${}^2$-completion of the space of $C^\infty$-sections
of the bundle $L$   and
$\gh$ be its (due to compactness of $M$) finite-dimensional
 closed subspace of holomorphic
sections.
Let $\ \Pi:\Ls\to\gh\ $ be the projection.
\begin{definition}
For $f\in\Cim$ the Toeplitz operator  $T_f$
is defined to be
$$ T_f:=\Pi\, (f\cdot):\quad\gh\to\gh\ .$$
\end{definition}
In words: We multiply the holomorphic section with the
differentiable function $f$. This yields only a differentiable section.
To obtain a holomorphic section again, we  project it back
to the subspace of global holomorphic sections.
{}From the point of view of Berezin's approach \cite{Berequ}, $T_f$
is the operator with contravariant symbol $f$.

The linear map
$$T:\Cim\to \End\big(\gh\big),\qquad  f\to T_f\ ,$$
is called the \BT\ quantization.
Recall that $(\Cim, \cdot,\{.,.\})$ is a Poisson algebra.
To define the Poisson bracket (i.e. a compatible Lie 
algebra structure) on $\Cim$ we use the K\"ahler form $\omega$
as symplectic form and define
$\{f,g\}:=\omega(X_f,X_g)$ where $X_f$ is the Hamiltonian
vector field assigned to $f\in\Cim$ given by $\omega(X_f,.)=df(.)$.
The \BT\ quantization map is
neither a Lie algebra homomorphism nor
an associative algebra homomorphism,
because in general
$$T_f\, T_g=\Pi\,(f\cdot)\,\Pi\,(g\cdot)\,\Pi\ne
\Pi\,(fg\cdot)\,\Pi\ .$$

Due to the compactness of $M$ this defines a map
from the commutative algebra of functions to a noncommutative
finite-dimensional (matrix) algebra.
A lot of information will get lost. To recover this
information one should consider not just the bundle $L$ alone but
all its tensor powers $L^m$ for $m\in\N_0$
and apply all the above constructions for every $m$.
In this way one obtains a family of
matrix algebras and maps
$$
\Tfm :\Cim\to \End\big(\ghm\big),\qquad  f\to \Tfm{}\ .$$
This infinite family should in some sense ``approximate'' the
algebra $\Cim$.(See \cite{BHSS}
for a definition  of such an approximation.)

If we group all $\Tfm$ together we obtain a map
\begin{gather}
\Cim\to \prod_{m\in\N_0}\End(\ghm)\subseteq 
\End(\prod_{m\in\N_0}\ghm),
\\
f\quad\mapsto\quad  T_f^{(*)}:=(\Tfm)_{m\in\N_0}\ .
\end{gather}
We will see later on that 
$\prod_{m\in\N}\ghm$ with a slight modification
(i.e. $\prod$ is replaced by $\bigoplus$) is 
the projective coordinate ring of the embedded $M$.
The operator $T_f^{(*)}$ is called the total Berezin-Toeplitz
operator. It operates on the projective coordinate ring.

It was shown by Bordemann, Meinrenken and Schlichenmaier
\cite{BMS} that this quantization scheme has the correct semi-classical
behavior and yields an associated star product (a deformation
quantization).
Denote by $||f||_\infty$ the sup-norm of $\ f\ $ on $M$ and by 
$||\Tfm||=\sup_{s\in\ghm, s\ne 0}\frac {||\Tfm s||}{||s||}$
 the operator norm on $\ghm$.
\begin{theorem}
\label{T:approx}
[Bordemann, Meinrenken, Schlichenmaier]\;\cite{BMS}

\noindent
(a) For every  $\ f\in \Cim\ $ there exists $C>0$ such that   
\begin{equation}
||f||_\infty -\frac Cm
\le||\Tfm||\le ||f||_\infty\ .
\end{equation}
In particular, $\lim_{m\to\infty}||\Tfm||= ||f||_\infty$.

\noindent(b) For every  $f,g\in \Cim\ $ 
\begin{equation}
\label{E:dirac}
||m\i[\Tfm,\Tfm[g]]-\Tfm[\{f,g\}]||\quad=\quad O(\frac 1m)
\ .
\end{equation}

\noindent(c) For every  $f,g\in \Cim\ $ 
\begin{equation}
||\Tfm \Tfm[g]-T^{(m)}_{f\cdot g}||\quad=\quad O(\frac 1m)\ .
\end{equation}
\end{theorem}
Let me recall the definition of a star product.
Let $\mathcal{A}=\Cim[[\nu]]$ be the algebra of formal power
 series in the
variable $\nu$ over the algebra $\Cim$. A product $\star$
 on $\mathcal {A}$ is 
called a (formal) star product if it is an
associative $\C[[\nu]]$-linear product such that
\begin{enumerate}
\item
\qquad $\mathcal{ A}/\nu\mathcal {A}\cong\Cim$, i.e.\quad $f\star g
 \bmod \nu=f\cdot g$,
\item
\qquad $\dfrac 1\nu(f\star g-g\star f)\bmod \nu = -\i \{f,g\}$,
\end{enumerate}
where $f,g\in\Cim$. 
We can also write 
\begin{equation}
\label{E:cif}
 f\star g=\sum\limits_{j=0}^\infty \nu^j C_j(f,g)\ ,
\end{equation}
with
$ C_j(f,g)\in\Cim$. The $C_j$ should be  $\C$-bilinear in $f$ and $g$.
The conditions 1. and 2.  can 
be reformulated as 
\begin{equation}
\label{E:cifa}
C_0(f,g)=f\cdot g,\qquad\text{and}\qquad
C_1(f,g)-C_1(g,f)=-\i \{f,g\}\ .
\end{equation}

\begin{theorem}
\label{T:star}
There exists a unique (formal) star product on $\Cim$
\begin{equation}
f \star g:=\sum_{j=0}^\infty \nu^j C_j(f,g),\quad C_j(f,g)\in
C^\infty(M),
\end{equation}
in such a way that for  $f,g\in\Cim$ and for every $N\in\N$  we have
with suitable constants $K_N(f,g)$ for all $m$
\begin{equation}
\label{E:sass}
||T_{f}^{(m)}T_{g}^{(m)}-\sum_{0\le j<N}\left(\frac 1m\right)^j
T_{C_j(f,g)}^{(m)}||\le K_N(f,g) \left(\frac 1m\right)^N\ .
\end{equation}
\end{theorem}
See \cite{Schlhab}, \cite{Schlbia95} and \cite{Schldef}.

It has a couple of nice properties, i.e.
(i) $1\star f=f\star 1=f$,
(ii) the selfadjointness $\overline{f\star g}=\overline{g}\star
\overline{f}$, and (iii) it admits a naturally defined trace
(see \cite{Schldef}).

As is shown in \cite{KarSchl}
the star product is a differential star product, i.e. the $C_j$ are
bidifferential operators and it has the property of 
``separation of variables" \cite{Karasep} (resp. it is of Wick type
\cite{BW}).
This says that it respects the holomorphic structure.
In more precise terms: if the star product is restricted to open
subsets the  star multiplication from the right 
with local holomorphic functions
is pointwise multiplication, and 
the  star multiplication from the left with local anti-holomorphic functions
is pointwise multiplication.

Let me close this section with two remarks.
\begin{remark}
More traditionally one considers the operator $Q$ of geometric 
quantization (with K\"ahler polarization) defined as
$Q=\Pi\circ P$ with 
$$
P:\Cim\to \End(\gul),\quad
f\mapsto P_f:=-\nabla_{X_f}+\i f\cdot id\ ,
$$
where $\gul$ is the space of $C^\infty$ sections of the
bundle $L$ and $\Pi$ is the projection onto the space of 
global holomorphic sections.
Now $Q_f\in\End(\gh)$.
Again one should consider $Q_f^{(m)}$ for all $m\in\N_0$.
For compact K\"ahler manifolds both quantization procedures are related via
the Tuynman relation. It reads  as
\begin{equation}
Q_f^{(m)}=\i\cdot T_{f-\frac 1{2m}\Delta f}^{(m)}=
\i\left(\Tfm-\frac 1{2m}T_{\Delta f}^{(m)}\right)\ .
\end{equation}
Hence, the $\Tfm$ and the $\Qfm$ have the same asymptotic 
behavior.
\end{remark}
\begin{remark}
There is another kind of embedding of the  manifold $M$ into 
projective space. It is the embedding using the coherent states
of Berezin-Rawnsley. This embedding turns out to be a special
case of the embedding considered 
at the beginning of this section where one uses a orthogonal 
basis of the sections, resp. (depending on the conventions) the conjugate
of it, see \cite{BerSchlcse} for details.
\end{remark}

\section{Projective varieties}\label{S:proj}
\subsection{The definition of a projective variety}\label{SS:proj}
Let $\K$ be an algebraically closed field and let us assume for 
simplicity that its characteristic is zero. 
Without any harm  the reader might even 
assume  $\K=\C$. 
The projective space $\P^n=\P^n(\K)$
is given as the space of lines through 
the origin in $\K^{n+1}$, i.e. 
as the equivalence classes of points in 
$\K^{n+1}\setminus\{0\}$ where two points $\alpha$ and  
$\beta$ are equivalent  if $\exists \lambda\in\K\setminus\{0\}$ with
$\beta=\lambda\cdot \alpha$.
The point $[\alpha]$ in projective space 
defined by the point 
$\alpha=(\alpha_0,\alpha_1,\cdots,\alpha_n)$, $\alpha\ne 0$ 
can be given by its (non-unique) homogeneous
coordinates 
$[\alpha]=(\alpha_0:\alpha_1:\cdots:\alpha_n)$.

Let $f\in\kno$ be a homogeneous polynomial of degree $k$.
As usual we obtain an associated $\K$-valued function on $\K^{n+1}$ 
by assigning to the point $\alpha\in\K^{n+1}$ the value
$f(\alpha)$ obtained by ``setting'' $X_i$ to be $\alpha_i$.
If $\beta=\lambda\alpha$ with $\lambda\in\K$, $\lambda\ne 0$,
is another point in the same
equivalence class as $\alpha$,  then we obtain 
$f(\lambda\alpha)=\lambda^kf(\alpha)$.
In particular, the induced function is only well-defined
on the whole projective space if $k=0$, i.e. $f$ is a 
constant.
But we also see that if $\alpha$ is a zero of $f$ then any
other element $\beta=\lambda\alpha$ will be a zero too.
Hence the zero-set
\begin{equation}
\Zero(f):=\{[\alpha]\in\P^n\mid f(\alpha)=0\}
\end{equation}
is a well-defined subset of $\P^n$.
The Zariski topology is the coarsest topology in which the sets $\Zero(f)$ are
closed subsets for all polynomials $f$, 
or equivalently for which the complements
$D_f=\P^n\setminus \Zero(f)$ are open sets.
Because the zero-sets of polynomials are also closed in the
``usual'' topology if the base  field is $\C$
the  sets which are closed  (open) in the Zariski topology are
closed(open) in the ``usual''  topology.
The Zariski topology has a number of quite unusual properties.
For example, it is not separated,
i.e. two distinct points do not necessarily have disjoint open 
neighborhoods. Even more is true: every non-empty Zariski open set $U$ is
automatically dense in $\P^n$.
\begin{definition}\label{D:proj}
(a) A subset $W$ of $\P^n$ is called a (projective) variety
if it is the set of common  zeros of finitely many homogeneous polynomials
$f_1,f_2,\ldots,f_m$ (which are not necessarily of the same degree)
\begin{equation}
W=\Zero(f_1,f_2,\ldots,f_m):=\{[\alpha]\in\P^n\mid f_i(\alpha)=0,\ i=1,\ldots,m\}.
\end{equation}

\noindent
(b) A variety is called a linear variety if it can be 
given as the zero-set of
linear polynomials.

\noindent
(c) A variety is called irreducible if every decomposition
\begin{equation}
W=V_1\cup V_2
\end{equation} 
with varieties $V_1$ and $V_2$  implies that
\begin{equation}
V_1\subseteq V_2\quad\text{ or }\quad
V_2\subseteq V_1\ .
\end{equation}
A variety which is not irreducible is called reducible.

\noindent
(d) A Zariski open set of a projective variety is called 
a quasiprojective (or sometimes just algebraic) variety.
\end{definition}
Note that 
some authors reserve the term variety for irreducible ones.
\begin{definition}\label{D:dimg}
Let $V$ be an irreducible variety, then its dimension 
$\dim V$ is
defined as the maximal length $n$ of chains of strict
subvarieties which are irreducible
\begin{equation}\label{E:chainvar}
\emptyset\ \subsetneqq \ V_0\  \subsetneqq\  V_1\ \subsetneqq\ \cdots\ 
\subsetneqq \ V_{n-1}\ \subsetneqq \ V_n=V\ .
\end{equation}
For arbitrary varieties the dimension is defined to be the
maximum of the dimensions of its irreducible subvarieties.
\end{definition}
Subvarieties of dimension 0
are called points, subvarieties of dimension 1 curves, etc.

Let $V$ be a projective variety 
i.e. $V=\Zero(f_1,f_2,\ldots,f_m)$.
Take $I=(f_1,f_2,\ldots f_m)$ to be the ideal generated 
by the polynomials $f_1,f_2,\ldots f_m$, i.e.
\begin{equation}\label{E:ideal}
I=\{\sum_{i=1}^m g_if_i\mid g_i\in\K[X_0,X_1,\ldots X_n],
i=1,\ldots,k\}\ .
\end{equation}
Obviously $V=\Zero(I)$.
Ideals which can be generated by homogeneous elements are 
called homogeneous ideals.
Hence, projective varieties can always be given as zero-sets of
homogeneous ideals. The converse is also true.
Clearly 
\begin{equation}
\Zero(I):=\{x\in\P^n\mid f(x)=0, \forall f\in I\}\ .
\end{equation}
is by the homogeneity of the generators a well-defined subset of
$\P^n$.
Because the polynomial ring is a Noetherian ring, i.e. every ideal 
$I$ 
can be generated (as ideal
in the sense of \refE{ideal}) by finitely many elements, e.g.
$I=(g_1,g_2,\ldots ,g_s)$, we get $\Zero(I)=\Zero(g_1,g_2,\ldots ,g_s)$
and hence $\Zero(I)$ is a projective variety in the sense
of \refD{proj}.
Any other set of generators of the ideal will define the same
zero-set. Even the ideal is not fixed uniquely by
$V$.
As a simple example one might consider 
the hyperplane $H=\Zero(I)$  
with $I=(X_0)$.
The same variety might be defined as 
 $H=\Zero(I'))$ with $I'=(X_0^2)$.
But note that $I'\subset I$.
One might expect that for a given
variety $V$ there is a largest ideal which still defines $V$.
This is indeed true.
\begin{definition}
Let $V$ be a projective variety, i.e. $V=\Zero(I)$ for some ideal $I$.
The vanishing ideal $\Ideal(V)$ is defined to be
\begin{equation}\label{E:vanish}
\Ideal(V):=\{f\in\Kon\mid f(x)=0,\forall x\in V\}\ .
\end{equation}
\end{definition}
The subset $\Ideal(V)$ is a homogeneous ideal and contains any
other defining ideal $I$ for $V$.
It can completely be described in algebraic terms.
For this we define for any ideal $I$ its radical ideal 
\begin{equation}\label{E:rad}
Rad(I):=\{f\in\Kon\mid \exists n\in\N:f^n\in I\}\ .
\end{equation}
If $I$ is homogeneous $Rad(I)$ will again be homogeneous.
We obtain
\begin{equation}
\Ideal(\Zero(I))=Rad(I)\ ,
\end{equation}
with only one exception in the case when $\Zero(I)=\emptyset$.
Note that $\emptyset$ corresponds to two homogeneous radical ideal,
the full ring $\Kon$ and the ideal $I_0:=(X_0,X_1,\ldots,X_n)$.
Note that the only possible zero of $I_0$ is the point
$0$ which is not an element of projective space.

There is another  warning necessary.
One might think that the dimension $r$ of a
variety is exactly $n-k$ if $k$ is the minimal number of necessary polynomials
to generate its  vanishing 
ideal $I$. Unfortunately this is not true. The only information one
has is that $r\ge n-k$, with equality if $k=1$. A variety is called a complete
intersection if indeed $r=n-k$.

Projective varieties are not always manifolds (of course not all
manifolds are projective varieties either).
Varieties have not necessarily to be smooth.
They might have singularities.
In \refS{sing} I will deal with singularities in more detail.
Here I would like to show some non-trivial examples of
singular varieties.
For this I  give a first  definition of a singular point.
Further definitions will follow in \refS{sing}.
\begin{definition}\label{D:sing}
Let $V=\Zero(f_1,f_2,\ldots,f_m)$ be a projective variety of dimension
$r$ in $\P^n$ 
with vanishing ideal $\Ideal(V)$ generated by 
the polynomials $f_1,f_2,\ldots ,f_m$.
Consider the $m\times(n+1)$-matrix (the Jacobi matrix)
\begin{equation}\label{E:jacmac}
J(X)=\left(
\begin{matrix}
\frac {\partial f_1}{\partial X_0}
& \frac {\partial f_1}{\partial X_1}
&\dots
& \frac {\partial f_1}{\partial X_n}
\\
\frac {\partial f_2}{\partial X_0}
& \frac {\partial f_2}{\partial X_1}
&\dots
& \frac {\partial f_2}{\partial X_n}
\\
\vdots&\vdots&\ddots&\vdots
\\
\frac {\partial f_m}{\partial X_0}
& \frac {\partial f_m}{\partial X_1}
&\dots
& \frac {\partial f_m}{\partial X_n}
\end{matrix}
\right)\ .
\end{equation}
A point $x$ on $ V$ with 
\begin{equation}
\rank J(x)<n-r
\end{equation} 
is called a singular point of the variety.
A point on a variety which is not a singular point is called a
regular point.
If $V$ has no singular points it is called a non-singular (or smooth, or
regular) 
variety.
\end{definition}
\begin{figure}
\begin{minipage}[b]{\linewidth}
\centering\epsfig{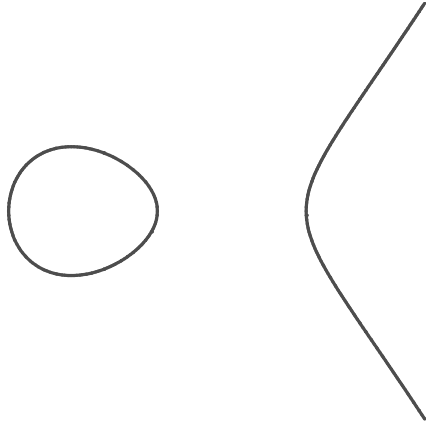}
\caption{A nondegenerate cubic curve}
\end{minipage}\hfill
\end{figure}
For a variety which is the union of two different 
subvarieties the points
where the subvarieties meet are always singular points. As a
typical example one might take the 
variety $V=\Zero(X_0X_1)$ in $\P^2$.
Then $V=\Zero(X_0)\cup\Zero(X_1)$ and the singular point is
the point $(0:0:1)=\Zero(X_0)\cap\Zero(X_1)$.
But even irreducible varieties can have singularities.
As an example let me consider the varieties $Y$ 
in $\P^2$ defined by 
irreducible cubic polynomials.
These polynomials can be written (after a suitable change of coordinates)
as 
\begin{equation}\label{E:cubic}
f(X,Y,Z)=Y^2Z-4X^3+g_2XZ^2+g_3Z^3\ ,
\end{equation} 
with certain elements $g_2,g_3\in\K$.
The variety 
$\Zero(f)$ is non-singular if and only if the coefficients
$g_2$ and $g_3$ are such that $g_2^3-27g_3^2\ne 0$.
One obtains in this way the {\it elliptic curves}.
See the Figure 1. Here the curve is defined over $\C$ and the
real-valued points are plotted.
The elliptic curves correspond to 1-dimensional complex tori,
see \refSS{embedd} below.
\begin{figure}
\begin{minipage}[b]{.44\linewidth}
\centering\epsfig{figure=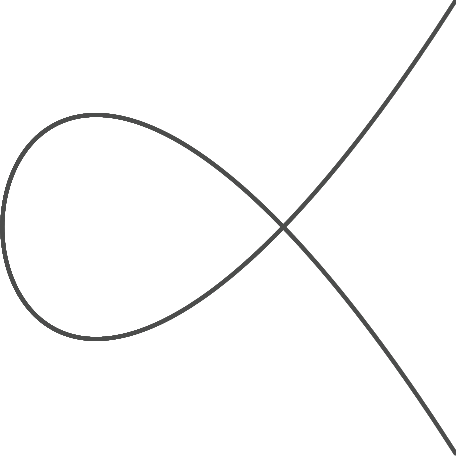,width=0.8\linewidth,clip=}
\caption{Nodal cubic}
\end{minipage}\hfill
\begin{minipage}[b]{.44\linewidth}
\centering\epsfig{figure=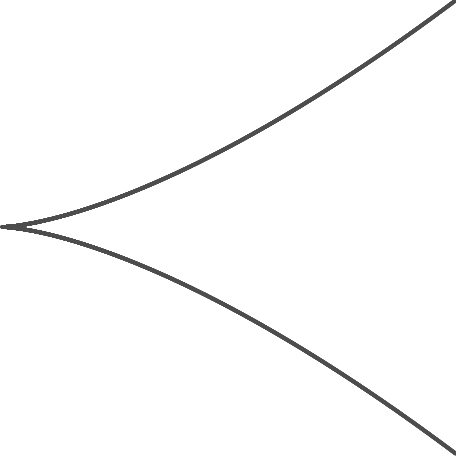,width=0.8\linewidth,clip=}
\caption{Cuspidal cubic}
\end{minipage}\hfill
\end{figure}
In the singular  cases we obtain two different types of curves.
The first one is the {\it nodal cubic}
$\Zero(Y^2Z-4X^2(X+Z))$, see Figure 2. 
The only singular  point is
the point $(0:0:1)$.
Moving along the curve we pass through the singular point twice, each time 
with a different tangent direction.
The second one is the {\it cuspidal cubic}
$\Zero(Y^2Z-4X^3)$, see Figure 3.
 Again the point $(0:0:1)$ is a singular point.
But now there is only one tangent direction at this point.

\subsection{Embeddings into Projective Space}\label{SS:embedd}
Let us now take the complex numbers $\C$ as base field $\K$.
With the complex topology $\P^n$ is a compact, $n$-dimensional
complex manifold. A coordinate covering is given 
by the ``affine'' sets 
\begin{equation}
D_{X_i}:=\Zero(X_i)=\{(\alpha_0:\alpha_1:\cdots:\alpha_n)\mid \alpha_i\ne 0\}
\ \cong\ \C^n
\end{equation}
for $i=0,1,\ldots,n$.
Every nonsingular projective variety is closed in the Zariski topology 
and hence also in the complex topology and hence a compact
submanifold of $\P^n$.
An abstract compact complex manifold $M$ is called 
a projective algebraic manifold  if there exists an injective holomorphic
embedding 
\begin{equation}
\Phi:M\to\P^n\ ,
\end{equation}
 such that $\Phi(M)\cong M$ as  complex manifolds.
The Theorem of Chow \cite[p.166]{GH} says that in this case
$\Phi(M)$ is a nonsingular projective variety, i.e. it can be given
as the zero-set of finitely many homogeneous polynomials.
This is even true in the strong sense that every meromorphic 
function on $\Phi(M)$ is a rational function
(i.e. it can be expressed as quotient of homogeneous polynomials
of the  same degree in $(n+1)$ variables),
every meromorphic differential is a rational differential, and every
holomorphic map between two embedded complex manifolds is an algebraic
map, i.e. can be given locally as a set of rational functions without
poles.

Let me illustrate this in the case of the above introduced elliptic curves.
Let $T=\C/\Gamma$ be the  
one-dimensional complex torus defined as the quotient of $\C$ by the
lattice 
$
\Gamma:=\{m+n\tau\in\C\mid m,n\in\Z\}$
for fixed $\tau\in\C$ with $\im\tau>0$.
The associated Weierstra\ss\ $\wp$-function and its derivative $\wp'$ are
doubly-periodic meromorphic functions with respect 
to the lattice $\Gamma$, i.e.
$$
\wp(z+\omega)=\wp(z),\quad
\wp'(z+\omega)=\wp'(z),\quad\text{for all}\ \omega\in\Gamma,
\ z\in\C\ .
$$
Hence they are meromorphic functions on  $T$.
The function $\wp$ fulfills
the famous differential equation
\begin{equation}
(\wp')^2=4\wp^3-g_2\wp-g_3
\end{equation}
with the Eisenstein series 
$$
g_2:=60\sum_{\omega\in\Gamma,\omega\ne 0} \frac 1{\omega^4},\qquad
g_3:=140\sum_{\omega\in\Gamma,\omega\ne 0} \frac 1{\omega^6}\ .
$$
An embedding of the torus into the projective plane is given by
\begin{equation}\label{E:wpdiff}
\Psi:T\to\P^2,\qquad
[z]\mapsto\begin{cases}
           (\wp(z):\wp'(z):1),&[z]\ne 0
            \\
           (0:1:0),&[z]= 0.
         \end{cases}         
       \end{equation}
Here $[z]=z\mod\Gamma$ denotes the point on the torus 
represented by $z\in\C$.
If one compares the differential equation \refE{wpdiff}
with the polynomial \refE{cubic} 
one sees that $\Psi(T)$ is a cubic
curve hence a projective variety
(indeed it is nonsingular).
Via $\Psi$ the meromorphic function $\wp$ corresponds to the rational
function $X/Z$ and $\wp'$ corresponds to $Y/Z$.
Note that the field of meromorphic functions on the torus 
consists of rational expressions in $\wp$ and $\wp'$.
For more details see, \cite{SchlRS}, p.34 and p.62 ff.

After this excursion let me return to the 
situation discussed in the \refS{quant}.
Let $M$ be a compact complex manifold and $\pi:L\to M$ a holomorphic line
bundle (not necessarily a quantum line bundle).
Choose a basis of the global holomorphic sections
$s_0,s_1,\ldots ,s_n\in\gh$.
For every point $x\in M$ there exists an open  neighborhood
$U$ of $x$ such that 
$L$ can be locally trivialized over $U$, i.e. that there is an
(holomorphic) bundle map
$\rho:L_{U}:=\pi^{-1}(U)\cong U\times\C$.
With respect to this trivialization the section $s_i$ can be given
by a local holomorphic function $\hat s_i:U\to\C$ defined by
$\rho(s_i(x))=(x,\hat s_i(x))$.
The map 
\begin{equation}
U\to\C^{n+1},\qquad y\mapsto \tilde\Phi(y)
:=(\hat s_0(y),\hat s_1(y),\dots,\hat s_n(y))
\end{equation}
is a holomorphic map. It depends not only on the basis chosen, but also
on the trivialization.
If $\rho'$ is a different trivialization defined over the open set $U$ 
 where  $\rho$ is defined (or a subset of it) 
then $\rho'\circ \rho^{-1}(x,\lambda)=(x,g(x)\lambda)$ 
with a holomorphic function 
$g:U\to\C$ nowhere vanishing on $U$.
The map  $\tilde \Phi':U\to\C^{n+1}$ corresponding to $\rho'$  
fulfills
$\tilde\Phi'(y)=g(y)\cdot \tilde\Phi(y)$.
Hence, 
$[\Phi(y)]:=\tilde\Phi(y)$ will be well-defined, i.e. not depend on the 
trivialization chosen if we assure that $\tilde\Phi(y)\ne 0$.
But $\tilde\Phi(y)=0$ if and only if
$s(y)=0$ for all sections $s\in\gh$.
Hence we obtain a well-defined
holomorphic map
\begin{equation}\label{E:lbmap}
\Phi:M\setminus\{y\in M\mid s(y)=0,\forall s\in\gh\}\ \xrightarrow{} \ 
\P^n\ ,
\end{equation}
obtained by glueing together the local maps $\tilde\Phi$.
A change of basis of
$\gh$ is given by an element of $\GL(n+1,\C)$.
The images of the two  mappings 
obtained by the two set of basis elements 
are related by the
corresponding $\PGL(n+1,\C)$ action.
Note that
if there exists a nontrivial section $s$ (i.e. $s\not\equiv 0$)
then the map $\refE{lbmap}$ 
is defined on a dense open subset of $M$. 
\begin{definition}\label{D:ample}
(a) A line bundle $L$ is called very ample if the map $\Phi$ 
(with respect to one and hence to all  set of basis elements) is an 
embedding.

\noindent
(b) A line bundle $L$ is called ample if there exists $m\in \N$ such that
$L^{\otimes m}$ is very ample.
\end{definition}
It follows that a compact complex manifold is projective algebraic if
it admits an ample line bundle.
The converse is also true. To see this we first study $\P^n$. Here we have
the tautological line bundle whose  fiber over the point
$[z ]$ is the complex line trough $0$ and the point $z\in\C^{n+1}$.
The hyperplane section bundle $H$ is the dual of the tautological line bundle.
Its space of global sections is generated by the coordinate functions
$X_0,X_1,\ldots, X_n$, i.e. it can be identified with the space of
linear polynomials in $(n+1)$ variables. All line bundles over $\P^n$ are
given 
as $H^m$ where this denotes for $m>0$ the $m$-th tensor power of $H$,
for $m<0$ the   $|m|$-th tensor power of the dual bundle of $H$, and
for $m=0$ the trivial bundle $\mathcal{O}$.
The space of 
global holomorphic sections of $H^m$ can canonically be identified 
with the space of homogeneous polynomials of degree $m$ in $(n+1)$ variables.
In particular, there exists no nontrivial sections for $m<0$. 
If $\Phi:M\to \P^n$ is a holomorphic map then the pullback 
$\Phi^*H$ is a holomorphic line bundle on $M$.
The space of global sections of $\Phi^*H$ is  generated by the pullback 
$\ \Phi^*(X_i)=X_i\circ \Phi\ $
of the global
sections $X_i$,  $i=0,1,\ldots, n+1$.
If $\Phi$ is a holomorphic embedding than $\Phi^*H$ is a very ample
line bundle. If the pull-backs of the $(n+1)$ sections $X_i$ stay linearly
independent then $\Phi$ is exactly given by the embedding defined via 
the bundle $\Phi^*H$. If not, then the  embedding defined via $\Phi^*H$
goes   into a  linear
subvariety of $\P^n$ of lower dimension, hence in a $\P^k$ for $k<n$.

Altogether
we see that the embeddings of $M$ into projective space  correspond to
very ample line bundles over $M$.
The pair $(M,L)$ where $M$ is a compact complex manifold and $L$ is a
very ample line bundle is called a {\it polarized projective 
algebraic manifold}.

Note that the same manifold considered with different $L$ may
``look'' quite differently. As a simple example take $M=\P^1$ and 
 $L=H$ then $\Phi: \P^1\to\P^1$ is the identity.
Now consider $L=H^2$, which gives an embedding into $\P^2$.
Let $X_0,X_1$ be the basis of the sections of $H$  then
$\ X_0^2,\ X_0X_1,\ X_1^2\ $ is a basis of $H^2$. If $(\alpha_0:\alpha_1)$ 
are homogeneous coordinates on $\P^1$ the image of $\P^1$ in
$\P^2$ is given  as
$$
\Phi(\P^1)=\{(\alpha_0^2:\alpha_0\alpha_1:\alpha_1^2)\in\P^2\mid
\alpha_0,\alpha_1\in\C\}=\Zero(X_1^2-X_0X_3)\ .
$$
The obtained subvariety is not linear anymore.
Nevertheless it is algebraically  isomorph to the linear variety $\P^1$.

Let us come back to the quantization condition.
Recall that the quantization  condition says that the Chern form of 
the quantum line bundle $L$
is essentially the K\"ahler form. But the K\"ahler form is a positive 
form, hence $L$ is a positive line bundle.
Kodaira´s embedding theorem says that a certain positive tensor power
of $L$ will give an embedding into projective space.
Hence  $L$ is an ample line bundle.
This implies that quantizable compact K\"ahler manifolds are always projective 
algebraic.
In \refSS{ring} we will see that the converse is also true.
\subsection{The projective coordinate ring}\label{SS:ring}
Let $V$ be a projective variety in $\P^n$ and $I=\Ideal(V)$ 
its vanishing ideal \refE{vanish}.
Recall that it is a homogeneous ideal.
\begin{definition}
The projective (or homogeneous) coordinate ring is the graded ring
\begin{equation}\label{E:ring}
\K[V]:=\K[X_0,X_1,\ldots X_n]/\Ideal(V) .
\end{equation}
\end{definition}
\noindent
For $V=\P^n$, we have $\Ideal(\P^n)=(0)$, hence
$$
\K[\P^n]:=\K[X_0,X_1,\ldots X_n]=\bigoplus_{m\ge 0}H^0(\P^n,H^m)
$$
is the full polynomial ring.

Inside $\K[V]$ the whole geometry of the variety $V$ is encoded.
For example the points correspond to maximal homogeneous ideals 
$M\subseteq I$ which are not identical to the ideal $(X_0,X_1,\ldots X_n)$.
Note that the only element
of $\K^{n+1 }$ which is a zero  of all $X_i$ is
$0$, which is not an element of projective space.
\begin{definition}\label{D:krull}
The Krull dimension $\dim R$ of a ring $R$ is defined
to be 
the maximal length $k$ of strict chains of prime ideals $P_i$
\begin{equation}
P_0\ \subsetneqq\   P_1\ \subsetneqq \ \cdots\  \subsetneqq
\  P_k\subseteqq\  R\ .
\end{equation} 
\end{definition}
Recall that an ideal $P$ is called a prime ideal if from $f\cdot g\in P$
 it follows that $f\in P$ or $g\in P$.
Clearly, for prime ideals $P$ we have $Rad(P)=P$, hence 
$\Ideal(\Zero(P))=P$.
Moreover, the variety $\Zero(P)$ is always irreducible.
Any chain \refE{chainvar} of irreducible subvarieties of an
irreducible variety
 gives  a chain of homogeneous prime ideals
\begin{equation}\label{E:chainideal}
\Kon\ \supsetneqq \ \Ideal(V_0)\  \supsetneqq\  \Ideal(V_1)
\supsetneqq \cdots\ \supsetneqq \ \Ideal(V_n)=\Ideal(V)\ .
\end{equation}
lying between the vanishing ideal of $V$ and the whole ring.
Passing to the quotient, i.e. to the coordinate ring
one obtains a chain of prime ideals of the coordinate ring $\K[V]$.
This works also in the opposite direction with the one exception
that to both the whole ring $\K[V]$ and to the ideal
 $(X_0,X_1,\ldots,X_n) \mod I(V)$ corresponds the empty set.
This implies
\footnote{
Note that for homogeneous coordinate rings 
to determine the Krull dimension it is enough to 
consider chains of homogeneous prime ideals.}
$$\dim V=\dim \K[V]-1\ .$$

Now let $\Phi:M\to\P^n$ be the embedding obtained 
via the quantum line bundle 
$L$, which we assume already to be very ample.
Let 
$I:=\Ideal(\Zero(\Phi(M))$
be the vanishing ideal of $\Phi(M)$.
We obtain $\Phi^*H=L$, $i^*X_i=s_i$ for $i=0,1,\ldots,n$  
for the sections $s_i$ used for the embedding,
and $\Phi^*(H^m)=(\Phi^*H)^m=L^m$.
In particular, the pull-backs of the global sections of $H^m$ generate the
space of global sections of  $L^m$.
But in general they will not be a basis. The algebraic 
relations between them are exactly given by the
elements of the ideal I.
The projective coordinate ring $\C[V]$ can be identified with 
$\bigoplus_{m\ge 0}H^0(M,L^m)$.

In \refS{quant} we have defined the Berezin-Toeplitz quantization 
map 
\begin{equation}
\label{E:BTmap}
\Cim\to \End\left(\prod_{m\in\N_0} H^0(M,L^m)\right),\quad
f\mapsto T_f^{(*)}=(T_f^{(m)})_{m\in\N_0}\ .
\end{equation}
Due to the fact that $T_f^{(*)}$ respects the grading
given by $m$, it can also
be considered as an element of 
$$
\End\left(\bigoplus_{m\in\N_0}H^0(M,L^m)\right)
$$ and is fixed by this restriction.
Hence, $T_f^{(*)}$ is an algebraic object operating on an algebraic
vector space which coincides with the coordinate ring.
The coordinate ring should be considered as the quantum Hilbert
space. Note that this set-up makes perfect sense also
for singular projective varieties.

Clearly, there is also a metric aspect in the theory. Our line  bundle comes
with a hermitian metric. On $\P^n$ we have 
the Fubini-Study K\"ahler form $\w_{FS}$ induced by the standard metric in 
$\C^{n+1}$. This defines a metric on the tautological bundle and 
by taking the inverse metric a hermitian metric $h_{FS}$
 on the hyperplane section bundle $H$. 
Suitable normalized it turns out that 
$H$ with $h_{FS}$ is the quantum line bundle of the K\"ahler manifold
$(\P^n,\w_{FS})$.
If $N$ is a closed submanifold of $\P^n$, i.e. a nonsingular projective 
variety and  $i:N\to\P^n$ is the embedding then  
the pair $(N,i^*\w_{FS})$ is a K\"ahler manifold with associated 
quantum line bundle $(i^*H,i^*h_{FS})$.
In particular, nonsingular projective varieties are always quantizable.
But note that if we start with a fixed K\"ahler manifold $(M,\w_M)$ with
very ample quantum line bundle $(L,h)$ and induced embedding
$\Phi:M\to \P^n$ then $(M,\Phi^*\w_{FS})$ is again a 
quantizable K\"ahler manifold with quantum line bundle
$(L\cong \Phi^*H,\Phi^*h_{FS})$. But in general 
we have for the two K\"ahler forms defined on the
same complex manifold $\Phi^*\w_{FS}\ne \w_M$.
We only know that they are cohomologous because they are representatives
of the Chern class of the same bundle $L$.
The question whether they coincide as forms has to do with the question whether
the embedding is a K\"ahler embedding. This is related to 
Calabi's diastatic function, respectively to Rawnsley's epsilon function.
I will not discuss this matter here, but see \cite{BerSchlcse} for 
a discussion and references to further results.

Via the metric the projective coordinate ring 
$\bigoplus_{m\ge 0}H^0(M,L^m)$ carries also a metric structure. 
To have a full description of the quantization also in the singular
case the metric should be studied in more detail.

\section{Singularities}\label{S:sing}
%
In the last section a point on a projective variety was 
called a singular point if the rank  of the matrix
\refE{jacmac} is less than expected (see \refD{sing}).
In this section I will give a different characterization
of singular points. In particular, it will turn out, that there
exist singularities which are better than others.

Clearly, the definition of a singular point as given in 
\refD{sing} is a local one. Hence it is enough to study the local situation.
For the local situation 
it is more convenient to consider affine varieties instead
of projective varieties.
If the projective space is 
replaced by an affine space the definitions work accordingly.
After choosing coordinates the $n$-dimensional affine space 
 is given as $\K^n$.
A subset $V$ of $\K^n$ 
is called an affine variety if there exists finitely many
polynomials $f_1,f_2,\ldots ,f_m\in\K[X_1,X_2,\ldots,X_n]$ such that
$$
V=\Zero(f_1,f_2,\ldots,f_m)
=\{\alpha\in\K^n\mid f_i(\alpha)=0,\ i=1,\ldots,m\}\ .
$$
If one
replaces homogeneous polynomials, homogeneous ideals, etc. by 
arbitrary polynomials,
arbitrary ideals, etc, the whole theory develops like in the projective case.
Again, let $I=(f_1,\ldots,f_m)$ be the ideal generated by the above
polynomials then $V=\Zero(I)$. 
Vice versa, given a variety $V$ in $\K^n$ we can define its vanishing ideal
\begin{equation}
\Ideal (V):=\{f\in \K[X_1,X_2,\ldots,X_n]\mid f(\alpha)=0,\forall \alpha\in V\}\ .
\end{equation}
With the same definition \refE{rad} of the radical ideal we obtain
$\Ideal(\Zero(I))=Rad(I)$.
The affine coordinate ring of the variety $V$ is defined to be
\begin{equation}
\K[V]:=\K[X_1,X_2,\ldots,X_n]/\Ideal(V)\ .
\end{equation}

The subset 
\begin{equation}
U^{(i)}:=\P^n\setminus\Zero(X_i)
=\{(\alpha_0:\alpha_1:\ldots:\alpha_n)\mid \alpha_i\ne 0\}
\end{equation}
of $\P^n$ is a Zariski open (and hence dense) subset of $\P^n$.
It can be identified with the affine space $\K^n$ via
the map
\begin{equation}\label{E:affi}
\Phi_i((\alpha_0:\alpha_1:\cdot:\alpha_n))\mapsto
(\frac {\alpha_0}{\alpha_i},\ldots,
\frac {\alpha_{i-1}}{\alpha_i},
\frac {\alpha_{i+1}}{\alpha_i},\ldots,
\frac {\alpha_n}{\alpha_i})\ .
\end{equation}
In this way $\P^n$ is covered by $(n+1)$ copies of affine $n$-space,
i.e. $\P^n=\cup_{i=0}^{n}U^{(i)}$.
Every projective variety can be covered by affine
varieties. 
Let $f_l(X_0,X_1,\ldots,X_n)$ for $l=1,\ldots,m$ be 
defining homogeneous
polynomials 
for the projective variety $V$.
Fix $i$ with $0\le i\le n$ and let $f_l^{(i)}$ 
be the polynomials in $n$ variables 
obtained from the $f_l$ by setting the variable $X_i$ to $1$.
Then $V^{(i)}=\Zero(f_1^{(i)},\ldots, f_m^{(i)})$
defines an affine variety. Via the map \refE{affi}
we can identify $V^{(i)}=V\cap U^{(i)}$.
Again $V=\cup_{i=0}^{n}V^{(i)}$.
In particular every point of the projective variety lies at least
in one of these affine varieties $V^{(i)}$.

In the following let $V$ be an affine variety. Again the dimension
$\dim V$ can be defined by \refD{dimg}.
This coincides with the Krull dimension
of the coordinate ring $\K[V]$, i.e. $\dim V=\dim\K[V]$.
In the affine case there is no subtraction of 1 necessary, because 
in this case 
there is a complete 1:1 correspondence between prime ideals and
irreducible subvarieties.
Note that if $Y$ is a irreducible projective variety all 
covering affine varieties $Y^{(i)}$ will be irreducible
affine varieties and vice versa. Additionally we have 
$\dim Y=\dim Y^{(i)}$ for non-empty  $Y^{(i)}$.

Singular points of affine varieties can be defined according to \refD{sing}
(of course now only $n$ variables will appear,
hence we get an $m\times n$ matrix) using  generators
$f_1,f_2,\ldots, f_m$ of the vanishing ideal of the variety.
If the affine variety $V^{(i)}$ comes from a projective variety $V$ as
described above then $x\in V$ corresponding to $\Phi(x)$ will be
a singular point of $V$ if and only if $\Phi(x)$ is  a singular point of
$V^{(i)}$.

There are some problems with  this definition of a singular point.
First, it is not a priori clear that it does not depend on the
chosen generators of the ideal $\Ideal(V)$.
Second, the starting point of the definition is a variety lying in some
affine space. But the singularity should be something intrinsic
to the variety and not depend on the affine space the variety is 
lying in.
It can be shown that indeed the notion does not depend on 
these choices. Nevertheless, a more intrinsic definition of
a singularity would be desirable.

There is such a  definition which deals with the local ring
$\O_{V,\alpha}$ 
of the point $\alpha$ on $V$.
This local ring is defined as follows.
Let $\alpha=(\alpha_1,\ldots,\alpha_n)\in\K^n$ be a point on $V$.
The vanishing  ideal of $\alpha$ in the polynomial
ring is the ideal
$M_\alpha=(X_1-\alpha_1,X_2-\alpha_2,\ldots,X_n-\alpha_n)$.
It is a maximal ideal. This says that every ideal which is strictly 
bigger than $M_\alpha$ is the whole polynomial ring.
The condition $\alpha\in V$ is equivalent to 
$M_\alpha\supseteq I=\Ideal(V)$.
{}From this it follows that $M_\alpha\bmod I$ is a maximal ideal of
$\K[V]$.
The local ring $\O_{V,\alpha}$ of the variety $V$ at the point $\alpha$
is defined as the {\it localization} of the 
ring $\K[V]$ with respect to the maximal ideal  $M_\alpha\bmod I$
(for simplicity we will denote it also by  $M_\alpha$)
\begin{equation}
\O_{V,\alpha}=\K[V]_{M_\alpha}\ .
\end{equation}
The localization is the ring of fractions where the denominators are
elements from the multiplicative set $\K[V]\setminus{M_\alpha}$.
It is a Noetherian local ring. Noetherian means that every ascending chain
of ideals becomes stationary (or terminates).
Local means that the ring has only one maximal ideal. Here the unique
maximal ideal is ${M_\alpha\bmod I}/{1}$
(again simply denoted by $M_\alpha$).
\begin{definition}\label{D:regular}
A local ring $R$ with maximal ideal $M$ 
is called a regular local ring if
\begin{equation}\label{E:regular}
\dim_{R/M}M/M^2=\dim R\ .
\end{equation}
\end{definition}
Due to the fact that $M$ is a maximal ideal of $R$ the 
quotient $R/M$ is a field, and $M/M^2$ is a vector space over 
$R/M$.
On the left hand side of \refE{regular} the vector space dimension
is meant, on the right hand side the Krull dimension is meant.

In our case where $R$ is the local ring coming from
the coordinate ring of a variety over an 
algebraically closed field $K$ we have $R/M\cong \K$.
\begin{definition}
A point $\alpha\in V$ is called a non-singular 
(or regular) point of $V$ if its local
ring $\O_{V,\alpha}$ is regular. If not it is called a singular point.
A variety is called { non-singular} (or smooth, or regular) 
if all points are regular. 
The subset of singular points of $V$ is denoted by
$Sing(V)$.
\end{definition}
$Sing(V)$ is always an algebraic subvariety of  codimension 
$\codim_V Sing(V)\ge 1$.
In particular the non-singular locus is a non-empty Zariski open subset of $V$.
If $V$ is irreducible then $Sing(V)$ is dense.
For irreducible $V$ the dimensions of all local rings $\O_{V,\alpha}$
are constant and equal to the dimension of $V$.

The $\K$-vector space $M_\alpha/M_\alpha^2$ is also called
the {\it Zariski cotangent space}, resp. its dual 
$(M_\alpha/M_\alpha^2)^*$ the {\it Zariski tangent space}.
In general
\begin{equation}
\dim_\K M_\alpha/M_\alpha^2\ge \dim \O_{V,\alpha}=dim V,
\end{equation}
where we assume for the last equality $V$ to be irreducible.
Hence we can also define the singular points to be the points
where the dimension of the Zariski tangent space is bigger than
the dimension of the (irreducible) variety.

Let me illustrate this in the case of cubic curves. In an
affine chart using the ideal 
\begin{equation}
I=\big(Y^2-4X(X-a)(X-b)\big)
\end{equation}
with $a,b\in \K$ 
we obtain the  cubic curves as $\Zero(I)$.
In this normalisation $\alpha=(0,0)$ lies on the cubic.
The cotangent space at $\alpha$ is given as
\begin{equation}\label{E:lquot}
M_\alpha/M_\alpha^2=(X,Y)\bmod I/ (X^2,Y^2,XY )\bmod I\ .
\end{equation}
{}From the relations given by $I$ we calculate
\begin{equation}
Y^2=4a b X-8(a+b) X^2+4X^3\bmod I\ .
\end{equation}
Hence
\begin{equation}
4a b X\in (X^2,Y^2,XY )\bmod I\ .
\end{equation}
If $a\cdot b\ne 0$ the element $Y$ will be enough to generate the
quotient \refE{lquot}. The tangent space will be one-dimensional
and $(0,0)$ will be a nonsingular point.
If either $a$ or $b$ equals $0$ the element $X$ will also be necessary to
generate the tangent space. Hence $(0,0)$ will be a singular point.

Given an arbitrary irreducible (projective or affine) 
variety $V$ then there exists 
a stratification of the singularity set $Sing(V)$ 
obtained in the following
manner.
Let $U=V\setminus Sing(V)$
be the Zariski open set of regular points then 
$V\setminus U$ is a closed subvariety. It can be decomposed  into 
finitely many irreducible varieties of dimension less than
$\dim V$
$$
V\setminus U=V_1^{(1)}\cup V_2^{(1)}\cup\cdots  \cup V_l^{(1)}\ .
$$
Again from this complement $Sing(V\setminus U)$  can be determined. 
It is a subvariety of the  variety $V$ from higher codimension.
This process can be repeated as long as there are singularities.
Because the  codimension strictly increases it has to stop after finitely many 
steps.

Guided by the algebraic properties of the local rings we have an 
hierarchy for  the types of singularities.
If $R$ is a ring without zero divisors then $Quot(R)$ 
is the ring whose elements are the fractions of elements in $R$
with denominator $\ne 0$.
\begin{definition}
A ring $R$ (without zero divisors) 
is called a normal ring if the  elements of
$Quot(R)$ which are solutions of algebraic 
equations with coefficients from $R$ and highest coefficient 1 lie
already in $R$.
\end{definition}
It is a classical result (Gau\ss\ Lemma) that $\Z$ is normal and also 
that  polynomial rings over fields are normal. 
\begin{definition}
Let $V$ be an irreducible variety.

\noindent
(a) A variety is normal at a point $\alpha\in V$  if  
the local ring $\O_{V,\alpha}$ is normal.

\noindent
(b) A variety is called normal if it is normal at every point 
$\alpha\in V$.

\noindent 
(c) A singular point is called a normal singular point, and the singularity
is called a normal singularity, if the variety is normal at this point.
\end{definition} 
A regular local ring is always normal. Hence, regular points are always
normal.
If $V$ is a normal variety it follows 
\begin{equation}
\codim_V Sing(V)\ge 2\ .
\end{equation}
This says that normal singular varieties are 
``less singular'' than 
generic singular varieties.
A lot of the singular varieties which appear as moduli spaces are normal.
Normal varieties behave from the point of functions defined on them
similar to nonsingular varieties.
For example, if $V$ is a variety of dimension $\ge 2$ and $x\in V$ 
a normal point then every regular function in $V\setminus \{x\}$ 
can be extended to a regular function on $V$.
Additionally normality is necessary to have a well-behaved theory of
(Weil-)divisors based on codimension 1 irreducible subvarieties.

For every irreducible affine variety $V$ with singularity set $Sing(V)$ 
there exists a normal  affine variety $\widetilde{V}$ and an
algebraic  morphism
$\pi:\widetilde{V}\to V$ such that 
\begin{equation}
\pi^{-1}(V\setminus Sing(V))\cong V\setminus Sing(V)\ .
\end{equation} 
The variety $\widetilde{V}$ is called the normalization of $V$.
It is obtained by a purely algebraic process, i.e. by taking 
the normal closure
of the coordinate ring in its quotient field.
This can be extended to the projective case too.

For algebraic curves normal points are always regular
(there is no space for codimension two subvarieties). Hence 
the normalization gives already a desingularization.
In the case of the above discussed singular cubic curves the normalization
is given by the line (affine, resp. projective).

The question arises whether it might be even possible to find for
every projective  variety $V$ with singularity set $Sing(V)$
 a nonsingular
projective variety $Y$ which coincides with $V$ outside $Sing(V)$,
and is minimal in a certain sense.
Such a $Y$ is called a desingularization and the whole process
is called a resolution of singularities. 
It was shown by Hironaka \cite{Hir} (see also \cite{Lip})
that there exists for projective varieties over 
fields of  characteristic zero (and this is the case we
are dealing with)
a resolution of singularities. 
More precisely, for every projective variety $V$ there exists a 
nonsingular projective variety $Y$ and a 
proper
\footnote{
the algebraic equivalent of a compact map}
algebraic map $f:Y\to V$ 
such that 
$f$ is an isomorphism over an open non-empty subset $U\subseteq V$, i.e.
$f^{-1}(U)\cong U$.

\section{Quotients}\label{S:quot}
In this section let us assume $\K=\C$ for the base field.
\subsection{Quotients in algebraic geometry}
Moduli spaces of geometric objects
are very often  varieties with singularities.
Typically, 
they are obtained starting from a smooth variety classifying 
the objects with respect to a certain ``presentation''.
To obtain the moduli space one has to ``divide out'' the different
presentations. 
Usually,  one has a group operating on the 
presentations and a candidate of the moduli space  
is given by  the quotient set under the group action, the orbit space.
Unfortunately, it is not always  possible to endow the 
quotient set with a compatible structure of a variety again.
Even if we allow the quotient to be an algebraic  scheme it will not
be possible.
In our context schemes will appear as ``varieties with multiplicities''.
It is quite reasonable that one should at least incorporate 
such objects in the theory. 
E.g. if we have two lines in the plane meeting at a point and we
move one line with the intersection point fixed, we
will nearly always have two lines. There is only one exception, 
if the moving line coincides with the fixed one.
In this case the configuration consists of one line.
But from the deformation point of view we should better count
this special  line twice, i.e. we should consider it as double line.
The language of schemes deals with such objects and even with much
more general ones. 
Nevertheless to avoid giving additional definitions 
I will still work on algebraic varieties (affine, projective, quasiprojective)
in the following.
But the reader should keep in mind that 
the language of algebraic schemes 
would be more appropriate for moduli problems.
See \cite{EHschemes} for an introduction to this field.

Let $X$ be an algebraic variety and   
$G$ a reductive algebraic group acting 
algebraically on $X$.
This means  that $G$ is the complexification of a maximal compact subgroup
$K$ of $G$. Of special importance (and this are the examples that the reader 
should keep in mind) are the groups $\GL(n)$, $\SL(n)$, and $\PGL(n)$.
As indicated above it is important to study
``quotients'' of $X$ under actions of the group $G$.
Mumford has given with his geometric invariant theory (GIT) \cite{GIT}
the principal tool to deal with such quotients.
\begin{definition}
A morphism of algebraic varieties $f:X\to Y$  
is called a good quotient  if
\begin{itemize}
\item[(1)]
$f$ is surjective and $G$-invariant, i.e. 
$f(gx)=f(x)$, for all $g\in G$ and $x\in X$,
\item[(2)]
$\left(f_*(\O_X)\right)^G=\O_Y$,
\item[(3)]
if $V$ is a $G$-invariant closed subset of $X$ then $f(V)$ is closed
in $Y$, and if $V_1$ and $V_2$ are $G$-invariant closed subsets of $X$ then
$$
V_1\cap V_2=\emptyset\implies f(V_1)\cap f(V_2)=\emptyset\ .
$$
\end{itemize}
\end{definition}
\noindent
In Condition (2) $\O_X$ and $\O_Y$ are the structure sheaves of the
varieties $X$ and $Y$. They are essentially nothing else as the 
sheaves of local regular  functions on $X$ and $Y$ respectively.
Condition (2) states that the local regular functions on $Y$ 
can be given as those local regular functions on $X$ 
which are constant along the fiber and  invariant under $G$.

A good quotient is a {\it categorical quotient}
in the sense that 
\begin{itemize}
\item[(1)]
 $f$ is constant on the orbits of the action,
\item[(2)]
for every algebraic variety $Z$ with a morphism $g:X\to Z$ which
is constant on the orbits of the $G$-action on $X$ 
there exist a unique morphism 
\newline
$\bar g:Y\to Z$ with
$g=\bar g\circ f$.
\end{itemize}

\begin{definition}
A morphism of algebraic varieties $f:X\to Y$  
is called a geometric quotient  if
\begin{itemize}
\item[(1)]
 $f:X\to Y$ is a good quotient,
\item[(2)]
for every $y\in Y$  the fiber $f^{-1}(y)$ consists exactly of one 
orbit under the group action.
\end{itemize}
\end{definition}
\noindent
If any of these quotients exists then they are unique.

For a good quotient there might exists fibers consisting of
several orbits under the group action and these orbits are
not necessarily closed (for a geometric quotient the orbits 
are always closed).

If we have a geometric quotient than the orbit space carries a structure
of an algebraic variety. But this condition is very often to strong to
be fulfilled. We have sometimes to assign several orbits to one
geometric point to obtain a geometric structure
and to end up (hopefully) with a good quotient. 
Mumford's concept of stability will help
to decide what to do.
Let $X\subseteq \P^n$ be a projective  algebraic variety and
$G$ a reductive algebraic group embedded into $\GL(n+1)$ with an 
action of $G$ on $X$ given by the standard linear action of $\GL(n+1)$ on
the points in $\P^n$.
\begin{definition}
(1) A point $x\in X$ is called semi-stable if and only if there exists
a non-constant $G$-invariant homogeneous polynomial 
$F\in \Con$ with $F(x)\ne 0$.

\noindent
(2) A point $x\in X$ is called stable if and only if
\begin{itemize}
\item[(a)] 
the dimension of the Orbit $O(x)$ under the $G$-action equals the
dimension of the group  and 
\item[(b)]
there exists a non-constant  $G$-invariant homogeneous polynomial 
$F\in \Con$ with $F(x)\ne 0$, and the action of $G$ on the zero set 
$X_F:=\{y\in X\mid F(y)=0\}$ is closed, i.e. if for
every $y_0\in X_F$ the orbit $O(y_0)$ is closed.
\end{itemize}
\end{definition} 
\noindent
The set of stable
\footnote{``Stable'' in the above introduced sense corresponds to 
``properly stable'' in the definition of Mumford.
Stability in his sense does not require the condition on
the dimension of the orbit.}
 points of $X$ under the above group action $G$ is
denoted by $X^s$, the set of semi-stable points is denoted by $X^{ss}$.
Both are open subset. Clearly,
$X^{s}\subseteq X^{ss}\subseteq X$.

Let me point out that the notion of stability might depend on the
embedding of the projective variety $X$ into projective space
and a corresponding linearization of the action of $G$.
Recall from \refSS{embedd} that for an abstract projective 
variety $X$ 
such an embedding is defined by the choice of a very ample line bundle
$L$ on $X$ and a choice of basis of its global sections.

\begin{theorem}\label{T:GIT}
Assume that $X^s$ is non-empty then there exists a projective algebraic 
variety $Y$  and a morphism $f_{ss}:X^{ss}\to Y$ such that
\begin{itemize}
\item[(1)]
$f_{ss}$ is a good quotient of $X^{ss}$ by $G$,
\item[(1)]
there exists an open subset $U\subseteq Y$  
such that $f_{ss}^{-1}(U)=X^s$ 
\newline
and $f_s:={f_{ss}}_{|X^s}:X^s\to U$ 
is a geometric quotient
of $X^s$ by $G$.
\end{itemize}
\end{theorem}
The good quotient is projective, but the geometric
quotient is as an open subset of something projective in 
general only quasi-projective.
If we interpret this in the opposite way, we see that we 
will need in general also 
non-stable (but still semi-stable) points to obtain 
projective (this means compact in the complex topology) moduli spaces.
Clearly, even if the projective variety
we started with was smooth 
there is no reason to expect that  the quotient will be smooth.

For more details one might consult \cite{GIT}, or for a more 
leisurely reading  \cite{New}.

\subsection{The relation  with the symplectic quotient}
In this subsection I want to quote results on the relation between the
quotients in algebraic geometry and the symplectic quotients.
The results are taken from Francis Kirwan's appendix
to the third  edition 
of Mumford's  book on GIT \cite{GIT}
and are due to Kirwan, Kempf and Ness.
More details and references can be found there.

Let $X$ be a nonsingular projective complex variety in $\P^n$,
and $G$ a reductive group acting linearly on $X$ via $\rho:G\to \GL(n+1)$.
If $K$ is any fixed maximal  compact subgroup of $G
$ then after a suitable choice of
coordinates the subgroup $K$ acts 
unitarily on $X$, i.e.  $\rho_{|K}:K\to \U(n+1)$.
Let $\mathfrak{k}$ be the Lie algebra of $K$,   $\mathfrak{k}^*$ 
its dual, and $\mu:X\to \mathfrak{k}^*$ the standard moment map
defined  for all $a\in\mathfrak k$ by
\begin{equation}
\mu(x)(a):=\frac {{}^t\overline{\hat x}\cdot \rho_*(a)\hat x}
                 {2\pi\i ||\hat x||}.
               \end{equation}
Here $\hat x\in\C^{n+1}\setminus\{0\}$ is any vector of
homogeneous coordinates representing the element $x\in X\subseteq \P^n$
and $\rho_*:\mathfrak k\to\mathfrak u(n+1)$  is the tangent map of 
$\rho_{|K}$.

In this situation the {\it symplectic quotient} 
(or {\it Marsden-Weinstein reduction}) is defined as $\mu^{-1}(0)/K$.
On the other hand we can define the good quotient 
(which is also a categorical quotient) of the 
semi-stable points $X^{ss}$ by $G$. By \refT{GIT}
it is a projective variety which is commonly also 
denoted by  $X//G$. It contains  
as open subset the geometric quotient $X^{s}/G$.
Immediately the following question arrises:
How are these quotients related?
\begin{theorem}\label{T:sympl}
(a) The point $x\in X$ is semi-stable if and only if 
$\,\overline{O_G(x)}\cap \mu^{-1}(0)\ne\emptyset$, 
 i.e. the closure of the orbit of $x$ 
under $G$ meets $\mu^{-1}(0)$.

\noindent
(b) $\mu^{-1}(0)\subseteq X^{ss}$.

\noindent
(c) The inclusion under (b) induces a homeomorphism
\begin{equation}
\mu^{-1}(0)/K\to X//G\ .
\end{equation}

\noindent
(d) If we denote by $\mu^{-1}(0)_{reg}$ the set of the $x\in\mu^{-1}(0)$
for which the tangent map  $d\mu_x$ of the moment map is 
surjective, then the homeomorphism under (c) restricts to
a  homeomorphism
\begin{equation}
\mu^{-1}(0)_{reg}/K\to X^s/G\ .
\end{equation}
\end{theorem}
Hence we see that as topological spaces the symplectic quotient is
isomorphic to the good quotient and the subspace of
the regular points 
$\mu^{-1}(0)_{reg}/K$ 
is isomorphic to the geometric quotient.

Things get  slightly more complicated if we consider also the
complex structure.
If $X$ is a compact K\"ahler manifold then the 
symplectic quotient $\mu^{-1}(0)_{reg}/K$  
carries a structure of a complex K\"ahler manifold away from the singularities.
But   $X^s/G$ carries also a complex structure from
the geometric quotient construction.
These two structures coincide on the subset where the symplectic
quotient has no singularities.
Hence if  $\mu^{-1}(0)_{reg}/K$ is a K\"ahler manifold the 
complex structure coincides with the complex structure of the geometric
quotient.
But there exists examples (e.g. given by Kirwan \cite[p.159]{GIT})
where the geometric quotient has no singularities but the symplectic
quotient has singularities in the sense that the K\"ahler structure
coming from the reduction process is singular at certain points.
Further examples of the relation between the structure 
of the singularities of the two type of quotients 
are given in the contribution of J. Huebschmann 
\cite{Huebsch} to this volume.

\providecommand{\bysame}{\leavevmode\hbox to3em{\hrulefill}\thinspace}


\end{document}